\magnification=1000

\font\titre=cmbx10 scaled\magstep1
\font\auteur=cmr10 scaled\magstep1

\def\square{{\vcenter{\hrule height.4pt
      \hbox{\vrule width.4pt height5pt \hskip5pt
           \vrule width.4pt}
      \hrule height.4pt}}}
\def\nl{\medskip\noindent}

\def\semidirect{\mathop
{\rlap{\kern 6.5pt\vrule height 5pt width .4pt depth 0pt} \times}}

\def\f{\varphi}
\def\Sym{{\rm Sym}}
\def\BN{{\bf N}}
\def\mod{{\rm mod}}
\def\Z{{\bf Z}}
\def\Im{{\rm Im}}
\def\Ker{{\rm Ker}}
\def\Conj{{\rm Conj}}
\def\e{\varepsilon}
\def\sg{{\rm sg}}
\def\order{{\rm order}}
\def\Perm{{\rm Sym}}
\def\Aut{{\rm Aut}}

%%%%%%%%%%%%%%%%%%%%%%%%%%%%%%%%%%%%%%%%%%%%%%%%%%%%%Coxeter graphs%%
\def\lijntje{\vrule height2.4pt depth-2pt width0.5in}

\def\vlijntje{\vrule height0.45in depth0.4pt width0.4pt}
\def\vlijn{\buildrel {\hbox to 0pt{\hss$\textstyle\circ$\hss}}\over\vlijntje}

\def\vtriple#1\over#2\over#3{\mathrel{\mathop{\kern0pt #2}\limits_{\hbox
to 0pt{\hss$#1$\hss}}^{\hbox to 0pt{\hss$#3$\hss}}}}
\def\rvtriple#1\over#2\over#3{\mathrel{\mathop{\kern0pt #2}\limits_{\hbox
to 0pt{\hss$#3$\hss}}^{\hbox to 0pt{\hss$#1$\hss}}}}

%diagrams
\def\vertex#1{\vtop{\hbox{$\circ$}\nointerlineskip\vskip 3pt\hbox{$\scriptstyle #1$}}}

\def\hedge#1{\hskip-.7pt\vrule height2.7pt width#1 depth-2.4pt\hskip-.6pt}

\def\Dn{\hbox{%
\vertex1\hedge{8ex}%
\vbox{%
  \hbox{%
    \vbox to 11pt{%
      \hbox{$\scriptstyle 2$}\nointerlineskip\vskip 3pt%
      \hbox{$\circ$}}}\nointerlineskip
  \hbox{\hskip.5ex\vrule height7ex width.3pt depth0.2ex}\nointerlineskip
  \hbox{\vertex3}}%
\hedge{8ex}\vertex4$\cdots\cdots$\vertex{n-1}\hskip -2.4ex\hedge{8ex}\vertex{n}}}

\def\An{\vtriple{\scriptstyle1}\over\circ\over{}\kern-1pt\lijntje\kern-1pt
\vtriple{\scriptstyle{2}}\over\circ\over{}\kern-1pt\lijntje\kern-1pt
\vtriple{\scriptstyle3}\over\circ\over{}
\cdots\cdots
\vtriple{\scriptstyle n-1}\over\circ\over{}\kern-1pt\lijntje\kern-1pt
\vtriple{\scriptstyle n}\over\circ\over{}\kern-1pt}

\def\Bn{\vtriple{\scriptstyle1}\over\circ\over{}
\kern-4pt{\buildrel {\scriptstyle 4}\over\lijntje}\kern-4pt\vtriple{\scriptstyle{2}}\over\circ\over{}
\kern-1pt\lijntje\kern-1pt
\vtriple{\scriptstyle 3}\over\circ\over{}
\cdots\cdots\vtriple{\scriptstyle n-1}\over\circ\over{}
\kern-1pt\lijntje\kern-1pt
\vtriple{\scriptstyle n}\over\circ\over{}\kern-1pt}

\def\Ffour{\vtriple{\scriptstyle1}\over\circ\over{}
\kern-1pt\lijntje\kern-1pt\vtriple{\scriptstyle{2}}\over\circ\over{}\kern-3pt{\buildrel 
{\scriptstyle 4}\over\lijntje}\kern-3pt\vtriple{\scriptstyle3}\over
\circ\over{}\kern-1pt\lijntje\kern-1pt
\vtriple{\scriptstyle4}\over\circ\over{}\kern-1pt}
\def\Eeight{\vtriple{\scriptstyle1}\over\circ\over{}\kern-1pt\lijntje\kern-1pt
\vtriple{\scriptstyle3}\over\circ\over{}\kern-1pt\lijntje\kern-1pt
\vtriple{\scriptstyle4}\over\circ\over{\buildrel
{\scriptstyle 2}\over\vlijn}\kern-1pt\lijntje\kern-1pt
\vtriple{\scriptstyle5}\over\circ\over{}\kern-1pt\lijntje\kern-1pt
\vtriple{\scriptstyle6}\over\circ\over{}\kern-1pt\lijntje\kern-1pt
\vtriple{\scriptstyle7}\over\circ\over{}\kern-1pt\lijntje\kern-1pt
\vtriple{\scriptstyle8}\over\circ\over{}\kern-1pt}
\def\Eseven{\vtriple{\scriptstyle1}\over\circ\over{}\kern-1pt\lijntje\kern-1pt
\vtriple{\scriptstyle3}\over\circ\over{}\kern-1pt\lijntje\kern-1pt
\vtriple{\scriptstyle4}\over\circ\over{\buildrel {\scriptstyle
2}\over\vlijn}\kern-1pt\lijntje\kern-1pt
\vtriple{\scriptstyle5}\over\circ\over{}\kern-1pt\lijntje\kern-1pt
\vtriple{\scriptstyle6}\over\circ\over{}\kern-1pt\lijntje\kern-1pt
\vtriple{\scriptstyle7}\over\circ\over{}}
\def\Esix{\vtriple{\scriptstyle1}\over\circ\over{}\kern-1pt\lijntje\kern-1pt
\vtriple{\scriptstyle3}\over\circ\over{}\kern-1pt\lijntje\kern-1pt
\vtriple{\scriptstyle4}\over\circ\over{\buildrel {\scriptstyle 2}
\over\vlijn}\kern-1pt\lijntje\kern-1pt
\vtriple{\scriptstyle5}\over\circ\over{}\kern-1pt\lijntje\kern-1pt
\vtriple{\scriptstyle6}\over\circ\over{}\kern-1pt}
\def\Dfive{\vtriple{\scriptstyle1}\over\circ\over{}\kern-1pt\lijntje\kern-1pt
\vtriple{\scriptstyle{2}}\over\circ\over{}\kern-1pt\lijntje\kern-1pt
\vtriple{\scriptstyle3}\over\circ\over{\buildrel
{\displaystyle \scriptstyle4}\over\vlijn}\kern-1pt\lijntje\kern-1pt
\vtriple{\scriptstyle5}\over\circ\over{}\kern-1pt}
\def\Dfour{
\vtriple{\scriptstyle1}\over\circ\over{}\kern-1pt\lijntje\kern-1pt
\vtriple{\scriptstyle{2}}\over\circ\over{\buildrel \scriptstyle3\over\vlijn}\kern-1pt\lijntje\kern-1pt
\vtriple{\scriptstyle4}\over\circ\over{}\kern-1pt}

%%%%%%%%%%%%%%%%%%%%%%%%%%%%%%%%%%%%%%%%%%%%%%%%%%%%%

\phantom{blabla}
\bigskip\bigskip\bigskip
\centerline{\titre On a theorem of Artin}
\bigskip
\centerline{\auteur Arjeh M. Cohen \& Luis Paris}
\bigskip
\centerline{\auteur April 2002}

%%%%%%%%%%%%%%%%%%%%%%%%%%%%%%%%%%%%%%%%%%%%%
\bigskip\bigskip
\centerline{\bf Abstract}

\bigskip
{\leftskip 1.5cm\rightskip 1.5 cm \noindent We determine the
epimorphisms $A \to W$ from the Artin group $A$ of type
$\Gamma$ onto the Coxeter group $W$ of type $\Gamma$, in case $\Gamma$ is
an irreducible Coxeter graph of spherical type, and we prove that the
kernel of the standard epimorphism is a characteristic subgroup of
$A$. This generalizes an over 50 years old result of Artin.
\par}

\nl
AMS Subject Classification: Primary 20F36; Secondary 20F55. Keywords:
Braid groups, Artin groups, Coxeter groups, characteristic subgroups, reflections.

%%%%%%%%%%%%%%%%%%%%%%%%%%%%%%%%%%%%%%%%%%%%%%%%%%%%
\bigskip\nl
{\titre 1. Introduction}

\bigskip
In 1947, Artin [Art] published the following two results on the braid
group $B_n$ on $n$ strings.

\nl
{\bf Theorem 1.1} (Artin, [Art]). {\it Suppose that $\f :B_n \to
\Sym_n$ is an epimorphism.  Then, up to an automorphism of $\Sym_n$,
the elements $\f (\sigma_i)$ $(i=1,\ldots,n-1)$
are the standard generators $(i,i+1)$ of
% the Coxeter group of type $A_{n-1}$, 
$\Sym_n$, with the following two exceptions for
$n=4$:
$$
\f (\sigma_1)=(1,2,3,4), \quad \f(\sigma_2)= (2,1,3,4), \quad \f(\sigma_3)= (1,2,3,4)
$$
and
$$
\f (\sigma_1)= (1,2,3,4), \quad \f (\sigma_2)= (2,1,3,4), \quad \f (\sigma_3)= (4,3,2,1),
$$
where $\sigma_i$ denotes the $i$-th standard generator of $B_n$.}

\nl
{\bf Theorem 1.2} (Artin, [Art]). {\it The pure braid group $PB_n$ is a characteristic subgroup 
of $B_n$.}

\bigskip
Let $\Gamma$ be a Coxeter graph with set of vertices $\{1, \dots, n\}$
(cf.\ [Bou] for terminology). For two objects $a,b$ and for $m \in
\BN$, define the word
$$
[a,b\rangle^m = \left\{\matrix{
&(ab)^{{m \over 2}}\hfill &{\rm if}\ m\ {\rm is\ even,}\hfill\cr
&(ab)^{{m-1 \over 2}}a\hfill &{\rm if}\ m\ {\rm is\ odd.}\hfill\cr}
\right.
$$
The {\it Artin group} of type $\Gamma$ is the group $A$ generated by $n$ elements $\sigma_1, \dots, 
\sigma_n$, and subject to the relations
$$
[\sigma_i, \sigma_j \rangle^{m_{i,j}} = [\sigma_j, \sigma_i \rangle^{m_{i,j}}
$$
for $1 \le i<j\le n$, and $m_{i,j} \neq \infty$, where
$(m_{i,j})_{i,j}$ is the Coxeter matrix of $\Gamma$. The Coxeter
system of type $\Gamma$ is denoted by $(W, S)$ with $S$ consisting of
the images $s_i$ of $\sigma_i$ $(i=1, \dots,n)$ under the natural
epimorphism $\mu: A \to W$.  The number of generators $n$ is called
the {\it rank} of the Artin group (and of the Coxeter group).  We say
that $A$ is {\it irreducible} if $\Gamma$ is connected, and that $A$
is of {\it spherical type} if $W$ is finite.  Define the {\it colored
Artin group} $CA$ of type $\Gamma$ to be the kernel of the standard
epimorphism $\mu: A \to W$.

Call two epimorphisms $\f_1, \f_2 : A \to W$ {\it equivalent} if there is an 
automorphism $\alpha$ of $W$ such that $\f_2= \alpha \circ\f_1 $. Members of the equivalence class 
of the standard epimorphism are called {\it ordinary} epimorphisms. The epimorphisms which are 
not ordinary are called {\it extraordinary}.

\bigskip
Note that the Artin group $A$ of spherical type $A_{n-1}$ is the braid
group $B_n$ of Theorem 1.1 and that the corresponding Coxeter group
$W$ coincides with $\Sym_{n}$.  In this case, $CA$ is the pure braid
group $PB_{n-1}$ of Theorem 1.2.  Thus, it makes sense to ask whether
Artin's results for type $A_{n-1}$ also hold for other Coxeter
graphs. In this paper, we prove the results corresponding to Theorems
1.1 and 1.2 for the irreducible Artin groups of spherical type
by proving the following two theorems.

\nl
{\bf Theorem 1.3.} {\it  Let $A$ be an irreducible Artin group of spherical 
type distinct from $A_3$ with corresponding Coxeter group $W$.
Then the equivalence classes of extraordinary epimorphisms $\f :A \to W$ 
are as given in Table 1.}
\bigskip
\halign{\qquad\qquad\qquad\hfil#\qquad\hfil&\hfil#\hfil\qquad&\hfil#\hfil\cr
Coxeter type\quad&epimorphism\quad& defined in Proposition\cr
$I_2(m)$ $(m\equiv0\pmod4)$&$\nu_m^{1},\nu_m^2$&2.2\cr
$B_n$ ($n$ even)&$\nu_n^1$&3.2\cr
$B_n$ ($n$ odd)&$\nu_n^1$, $\nu_n^2$&3.2\cr
$D_n$ ($n$ odd)&$\nu_n^1$, $\nu_n^2$&4.2\cr
$H_3$ &$\nu_3^1$&5.2\cr
}
\medskip 
\centerline{\kern-100pt Table 1. 
The extraordinary epimorphisms up to equivalence.}

\nl
{\bf Theorem 1.4.} {\it Let $A$ be an irreducible Artin group of
spherical type. Then $CA$ is a characteristic subgroup of $A$.}

\bigskip
Each `family' of irreducible Artin groups of spherical type is treated
separately.  Our strategy for the proof of Theorem 1.3 is to analyze
the epimorphisms $A\to W$ in much the same way Artin did for type
$A_n$, and to study epimorphisms $A\to \Sym_n$ for types $B_n$,
$D_n$.  In order to prove Theorem 1.4, we establish that the kernel of
every surjective composition $\mu\circ\Phi$ of $\mu$ with an
automorphism $\Phi$ of $A$ coincides with $CA$. For this purpose, we study an
action of $A$ on the set $\Z\times T$, where $T$ is the set of reflections
of $W$, which is compatible (under $\mu$) with the action of $W$ on
the set $\{\pm1\}\times T$ of abstract walls of $(W,S)$ (Lemma 3.7).

The groups of rank two are treated in Section 2, Artin groups of type
$B_n$ in Section 3, Artin groups of type $D_n$ in Section 4, and Artin
groups of type $H_3$, $H_4$, $F_4$, $E_6$, $E_7$, $E_8$ in Section
5. As noted above, type $A_n$ is dealt with in [Art].  The
classification of the epimorphisms $A \to W$ for $A$ of type $B_n$ and
$D_n$ is almost done in [Lin]; our contribution consists of Lemmas 3.5
and 4.4 and leads to a complete classification.  Our approach to these
two families is the same as Lin's, except that we do not use any
homological argument.

\nl
{\bf Acknowledgments.} The second author would like to thank Jean
Michel who introduced him to GAP and to the package `Chevie' of GAP.
%%%%%%%%%%%%%%%%%%%%%%%%%%%%%%%%

\bigskip\nl
{\titre 2. Artin groups of rank two}

\bigskip
We start this section with a lemma which will be used throughout the whole paper. Its proof is direct 
from the definition of a Coxeter group.

\nl
{\bf Lemma 2.1.} {\it Let $A$ be an Artin group of spherical type. If
all the images of the standard generators under an epimorphism $\f: A
\to W$ are involutions, then $\f$ is ordinary.} $\ \square$

\bigskip
Recall from [Bou]
that $I_2(m)$ denotes the Coxeter graph on two nodes which are joined by an $m$-valued edge.
It is irreducible whenever $m\ge3$.
By abuse of notation, we will use $I_2(3)$, $I_2(4)$, and $I_2(6)$ to denote the Coxeter graph of 
type $A_2$, $B_2$, and $G_2$, respectively.

\nl
{\bf Proposition 2.2.} {\it Let $m\ge 3$. For the Coxeter graph
$I_2(m)$, each epimorphism $A \to W$ is ordinary if $m \not\equiv 0
(\mod\ 4)$. If $m \equiv 0 (\mod\ 4)$, then there are two equivalence
classes of extraordinary epimorphisms. These are represented by
$\nu_m^1,\nu_m^2: A \to W$, where}
\bigskip\centerline{\vbox{\halign{\hfill$#$&\ $#$\ &$#$\hfill\quad&\hfill$#$&\ $#$\ &$#$\hfill\cr
\nu_m^1(\sigma_1)&=&s_1,&\nu_m^1(\sigma_2)&=&s_2s_1,\cr
\noalign{\smallskip}
\nu_m^2(\sigma_1)&=&s_1s_2,&\nu_m^2(\sigma_2)&=&s_2.\cr}}}
\bigskip

\noindent
{\bf Proof.} Let $r=s_2s_1$. Let $C$ be the subgroup of $W$ generated
by $r$. It is a cyclic normal subgroup of $W$ of order $m$, and
$W=C\semidirect \langle\sigma_1\rangle\cong
\Z/m\Z\semidirect\Z/2\Z$.  Let $\f: A \to W$ be an epimorphism, and
let $t_1=\f(\sigma_1)$ and $t_2=\f (\sigma_2)$. Four cases have to be
considered.

\smallskip\noindent
{\it Case 1:} $t_1, t_2 \in s_1 C$. Then $t_1$ and $t_2$ are involutions, and $\f$ is ordinary by 
Lemma 2.1.

\smallskip\noindent
{\it Case 2:} $t_1, t_2 \in C$. Then $\Im(\f) \subset C$, and $\f$ is not surjective. So, this 
case does not occur.

\smallskip\noindent
{\it Case 3:} $t_1 \in s_1C$ and $t_2 \in C$. If $m$ is odd, then
$t_1$ and $t_2$ are conjugate (since $\sigma_1$ and $\sigma_2$ are
conjugate), a contradiction with $t_2\in C$ and $t_1\not\in C$ as $C$
is a normal subgroup of $W$. So, $m$ is even. The fact that $t_1$ and
$t_2$ generate $W$ implies that $t_2$ is of order $m$, so, up to an
automorphism of $W$, we can assume that $t_1=s_1$ and
$t_2=r$. Consider the equality
$$
t_1t_2t_1t_2=s_1rs_1r=1.
$$
If $m \equiv 2(\mod\ 4)$, then
$[t_1,t_2\rangle^m=[t_2,t_1\rangle^m$ implies $s_1r=rs_1$, but this
last equality does not hold. Hence, $m \equiv 0 (\mod\ 4)$. Now $\f$
is indeed a homomorphism as
$$
[t_1, t_2 \rangle^m= [t_2,t_1\rangle^m=1,
$$
and $\f$ is equivalent to $\nu_m^1$.

\smallskip\noindent
{\it Case 4:} $t_1 \in C$ and $t_2 \in s_1C$. 
Then Case 3 applies with the indices $1$ and $2$
interchanged. Therefore $m \equiv 0 (\mod\ 4)$, and $\f$ is 
equivalent to $\nu_m^2$. $\ \square$

\nl
Recall that $\mu: A \to W$ denotes the standard 
epimorphism mapping $\sigma_i$ to $s_i$.
\nl
{\bf Lemma 2.3.} {\it If $\Phi \in \Aut(A)$ is such that $\mu\circ \Phi$ is ordinary,
then $\Phi$ preserves $CA$.}

\nl
{\bf Proof.} As $\mu\circ\Phi: A \to W$ is ordinary, 
there is an automorphism $\alpha$ of $W$
such that $\mu\circ\Phi = \alpha\circ\mu$. Now
$$
\Phi^{-1}(CA)=\Phi^{-1}( \Ker (\mu))= \Ker(\mu\circ\Phi)=\Ker (\alpha\circ\mu)=\Ker(\mu)= CA.
\quad \square$$

\nl
{\bf Proposition 2.4.} {\it Let $m \ge 3$. For the Coxeter graph $I_2(m)$, the colored Artin 
group $CA$ is a characteristic subgroup of $A$.}

\nl
{\bf Proof.} Let $\Phi: A \to A$ be an automorphism.  In view of Lemma 2.3,
we may assume that $\mu \circ \Phi$ is an extraordinary epimorphism. In particular, by
Proposition 2.2, $m \equiv 0 
(\mod\ 4)$. Let $\Delta= [\sigma_1, \sigma_2 \rangle^m \in A$, and let $\delta= [s_1,s_2 
\rangle^m \in W$. Since $m$ is even, by [BS], the center of $A$ is the infinite cyclic subgroup 
generated by $\Delta$, so $\Phi( \Delta) = \Delta^\varepsilon$, where $\varepsilon \in \{ \pm 1 
\}$, whence
$$
(\mu \circ \Phi) (\Delta) = \mu( \Delta^\varepsilon) = \delta^\varepsilon = \delta .
$$
On the other hand, by Proposition 2.2, there exist an automorphism $\alpha$ of $W$ and $i \in 
\{1, 2\}$ such that $\mu \circ \Phi = \alpha \circ \nu_m^i$. But
$$
\nu_m^1(\Delta)=[s_1,s_2s_1 \rangle^m =1\quad\hbox{ and }\quad
\nu_m^2(\Delta)=[s_1s_2,s_2 \rangle^m=1,
$$
whence $ \delta = \mu\circ \Phi(\Delta) = (\alpha \circ \nu_m^i)
(\Delta)=1$. This contradicts $\delta\ne1$, so we are done. $\ \square$

%%%%%%%%%%%%%%%%%%%%%%%%%%%%%%%%%
\bigskip\nl
{\titre 3. Artin groups of type $B_n$}

\bigskip
For a group $G$ and an element $g \in G$, we denote by $\Conj_g: G \to G$ the inner automorphism 
$h \mapsto ghg^{-1}$. Call two epimorphisms $\f_1, \f_2 : G \to H$ {\it conjugate} if there is $h 
\in H$ such that $\f_2= \Conj_h \circ \f_1$.

We shall label the nodes of the Coxeter graph $B_n$ as follows:
$$\Bn$$

Let $A$ be an Artin group of type $B_n$, $n \ge 3$, and let $W$ be its corresponding Coxeter group. 
The first step in our investigation of the epimorphisms $A \to W$ is to classify the epimorphisms 
$A \to \Sym_n$ up to conjugation. This has also been
carried out by Zinde [Zin], and the answer is as 
follows.
Define the standard epimorphism $\eta: A \to \Sym_n$ by
$$
\eta(\sigma_1)=1, \quad \eta (\sigma_i) = (i-1,i)\quad{\rm for}\ i=2, \dots, n.
$$

\noindent
{\bf Proposition 3.1} (Zinde, [Zin]). {\it Let $n \ge 3$. For the Coxeter graph $B_n$, every 
epimorphism $A \to \Sym_n$ is conjugate to $\eta$ if $n \neq 4$ and $n \neq 6$.

If $n=4$, then there are precisely three conjugacy classes of
epimorphisms; these are represented by $\eta, \zeta_4^1,
\zeta_4^2: A \to \Sym_4$, where
$$
\zeta_4^1(\sigma_1)= 1, \quad \zeta_4^1(\sigma_2)= (1,2,3,4), \quad \zeta_4^1(\sigma_3)= 
(2,1,3,4), \quad \zeta_4^1(\sigma_4)= (1,2,3,4),
$$
and
$$
\zeta_4^2(\sigma_1)= 1, \quad \zeta_4^2(\sigma_2)= (1,2,3,4), \quad \zeta_4^2(\sigma_3)= 
(2,1,3,4), \quad \zeta_4^2(\sigma_4)= (4,3,2,1).
$$

If $n=6$, then there are precisely two conjugacy classes of
epimorphisms; these are represented by $\eta, \zeta_6: A \to
\Sym_6$, where
\bigskip
\centerline{\vbox{\halign{\hfill$#$&\ $#$\ &$#$\hfill\quad&\hfill$#$&\ $#$\ 
&$#$\hfill\quad&\hfill$#$&\ $#$\ &$#$\hfill\cr
\zeta_6(\sigma_1)&=& 1,& \zeta_6(\sigma_2) &=& (1,2)(3,4)(5,6),& \zeta_6(\sigma_3) &=& 
(1,5)(2,3)(4,6),\cr
\noalign{\smallskip}
\zeta_6(\sigma_4) &=& (1,3)(2,4)(5,6),& \zeta_6(\sigma_5) &=& (1,2)(3,5)(4,6), &\zeta_6(\sigma_6) 
&=& (1,4)(2,3)(5,6). \quad\cr}}}

\noindent
Moreover, if we set

\centerline{\vbox{\halign{\hfill$#$&\ $#$\ &$#$\hfill\quad&\hfill$#$&\ $#$\ 
&$#$\hfill\quad&\hfill$#$&\ $#$\ &$#$\hfill\cr
r_1&=&1, &r_2&=&s_2s_4s_6, &r_3&=&s_2s_3s_2s_5s_4s_3s_2s_6s_5,\cr
r_4&=&s_3s_2s_4s_3s_6, &r_5&=&s_2s_5s_4s_6s_5, &r_6&=&s_2s_3s_2s_4s_3s_2s_6,\cr}}}
then $r_i=\zeta_6(\sigma_i)$ 
for all $i=1, \dots, 6$, and the automorphism $\alpha$ of $\Sym_6$ determined by
$(i-1,i)\mapsto r_i$ $(i=2,\ldots,6$) 
satisfies $\zeta_6=\alpha\circ\eta$.
$\ \square$}
\bigskip

The proof of Proposition 3.1 can be recovered from Artin's proof of 
[Art, Thm.1] by an easy adaptation.

\bigskip
We need a more precise description of the well-known isomorphism $W\cong
(\Z/2\Z)^n\semidirect\Sym_n$.  Put $c_1=s_1$ and, for $i=2,\ldots,n$,
denote by $c_i$ the element ${s_i\cdots s_2}s_1{s_2\cdots s_i}$ of
$W$.  Then $C=\langle c_1,\ldots,c_n\rangle\cong (\Z/2\Z)^n$ is the
kernel of the natural epimorphism $p: W = C \semidirect S \to \Sym_n$,
where $S = \langle s_2,\ldots,s_n\rangle\cong \Sym_n$.  It is
convenient to identify $c_1^{a_1}\cdots c_n^{a_n}\in C$ with
$[a_1,\ldots,a_n]\in (\Z/2\Z)^n$, so that each element of $W$ can be
viewed as the product of an element from $(\Z/2\Z)^n$
and a permutation in $\Sym_n$. 
Observe that the longest element $\delta = (s_1\cdots s_n)^n$ of $W$ coincides with $c_1\cdots c_n$.

% We use Proposition 3.1 to prove the
%following.

\nl
{\bf Proposition 3.2.} 
{\it Let $n\ge3$. For the Coxeter graph $B_n$, 
representatives of the
equivalence classes of 
extraordinary epimorphisms $A \to W$
are
\item{(i)} $\nu_n^1$,
if $n$ is even,
\item{(ii)}
$\nu_n^1$ and $\nu_n^2$,
if $n$ is odd,
\nl
where
$\nu_n^1$ is the homomorphism
$A \to W$ determined by
$$\eqalign{
\nu_n^1(\sigma_1)&=c_1\cr
\nu_n^1(\sigma_i)&= c_{i-1}s_i \quad {\it for}\ i=2, \dots n,\cr}
$$
and, for $n$ odd,
$\nu_n^2$ is the homomorphism $A\to W$ determined by
$$\eqalign{
\nu_n^2(\sigma_1) &= \delta\cr
\nu_n^2(\sigma_i) &= c_{i-1}s_i \quad {\it for}\ i=2, \dots, n.\cr}
$$

}

\bigskip
Before proving
Proposition 3.2, we establish three lemmas.
%3.3, 3.4, and 3.5 
Observe that $\eta = p\circ \mu$.

\nl
{\bf Lemma 3.3.} {\it Let $\f: A \to W$ be an extraordinary epimorphism such that $p \circ \f$ is 
the standard epimorphism $\eta: A \to \Sym_n$. Then, up to an automorphism of $W$, we have}
$
\f (\sigma_i) =c_{i-1}\cdot s_i \quad {\it for}\ i=2, \dots, n
$.

\noindent
{\bf Proof.} Let $t_i= \f (\sigma_i)$, $i=1,2, \dots, n$. Note that
$t_1 \in C$ in view of the assumption that $p\circ\f =
\eta$. Thus, either $t_1=1$ or $t_1$ is an involution. Note also that
$t_i$ is conjugate to $t_2$ for $i=3, \dots, n$.

Assume first that $t_2$ is an involution. For $i=3, \dots, n$, the element $t_i$,
being conjugate to $t_2$, is also an involution. Furthermore,
$t_1$ cannot be the identity, for otherwise we would have an
epimorphism $\Sym_n \to W$, $s_i \mapsto t_{i+1}$ ($i=1, \dots,
n-1$). Hence, $t_1$ is also an involution, and, by Lemma 2.1, $\f$ is
ordinary. This is contrary to the hypotheses, so $t_2$ is not an involution.

As $p\circ\f=\eta$, for $i=2, \dots, n$ we can write $t_i=v_i\cdot
s_i$ with $v_i=c_1^{a_{i,1}} \cdots c_n^{a_{i,n}} \in C$.  The
inequalities $t_i^2 \neq 1$ imply $a_{i,i} \equiv 1+a_{i,i-1}\pmod{2}$
for $i=2, \dots, n$. The equalities $t_i t_2=t_2 t_i$ imply $a_{2,i}
\equiv a_{2,3}\pmod2$ for $i=4, \dots, n$.  Finally, the equalities
$t_i t_{i+1} t_i = t_{i+1} t_i t_{i+1}$ ($i=2, \dots, n-1$) imply
$a_{i,j} \equiv a_{2,3}\pmod2$ for $i=2, \dots, n$ and $j \neq
i-1,i$. Write $a = a_{2,3}$ and $a_i= a_{i,i-1}$ for $i=2, \dots,
n$. Then $$v_i= c_1^a \cdots c_{i-2}^a c_{i-1}^{a_i} c_i^{1+a_i}
c_{i+1}^{a} \cdots c_{n}^{a}.$$

Suppose $a\equiv1\pmod2$. Let $\alpha_0$ be the automorphism of $W$ determined by
$$\eqalign{
\alpha_0(s_1) &=s_1,\cr
\alpha_0(s_i) &= \delta s_i \quad{\rm for}\ i=2, \dots, n.\cr}
$$
Then
$$
(\alpha_0 \circ \f)(\sigma_i) = \alpha_0(t_i) = c_{i-1}^{1+a_i}c_{i}^{a_i} s_i
\quad {\rm for}\ i=2, \dots, n,
$$
so we can assume $a=0$.

Suppose $a_n\equiv \cdots \equiv a_{i+1}\equiv1\pmod2$ and $a_i\equiv0\pmod2$. Then
$$\eqalign{
(\Conj_{c_{i-1}} \circ \f) (\sigma_j) = \Conj_{c_{i-1}}(t_j) &= t_j \quad {\rm if}\ j \neq i-1,i,\cr
(\Conj_{c_{i-1}} \circ \f) (\sigma_j) = \Conj_{c_{i-1}}(t_j) &= c_{j-1} c_{j}  
t_j \quad {\rm if}\ j = i-1,i,\cr}
$$
so we can assume $a_n\equiv \cdots \equiv a_{i+1}\equiv a_i\equiv1\pmod2$. 
An iteration of this argument shows that we can assume $a_i\equiv1\pmod2$ for
$i=2, \dots, n$. But then
$v_i=c_{i-1}$. % up to an automorphism of $W$. 
$\ \square$

\nl
{\bf Lemma 3.4.} {\it Let $\f: A \to W$ be an extraordinary
epimorphism such that $p \circ \f=\eta$.
\item{(i)}
If $n$ is even, then $\f$ is equivalent to
$ \nu_n^1$.
%$$\eqalign{
%\f (\sigma_1) &= c_1,\cr
%\f (\sigma_i) &= c_{i-1} \cdot s_i \quad{\it for}\ i=2, \dots n.\cr}
%$$
\item{(ii)}
If $n$ is odd, then $\f$ is equivalent to
$ \nu_n^1$ or to $ \nu_n^2$.
%$$\displaylines{
%{\it either}\quad \f (\sigma_1)= c_1 \quad {\it or} \quad \f (\sigma_1) = \delta,\cr
%{\it and} \quad \f (\sigma_i)= c_{i-1} \cdot s_i\quad {\it for}\ i=2, \dots, n.\cr}
%$$
%
}

\nl
{\bf Proof.} Let $t_i= \f(\sigma_i)$, for $i=1,2, \dots, n$. By Lemma 3.3, we can assume $t_i = 
c_{i-1}\cdot s_i$ for $i=2, \dots, n$. Write $t_1= c_1^{a_1}c_2^{a_2} \cdots c_n^{a_n} \in C$.
The equalities $t_i t_1=t_1 t_i$ ($i=3, \dots, n$) imply $a_i\equiv a_2\pmod2$ for $i=3, \dots, n$, so, 
$t_1= c_1^{a_1}c_2^{a_2} \cdots c_n^{a_2} $.

Consider the epimorphism $\sg: W \to \{\pm1\}$, $s_i \mapsto -1$. 
It satisfies $\sg(t_1) \equiv (-1)^{a_1+(n-1)a_2}$ and
$\sg(t_i) =1$ for $i=2,\ldots,n$.
So, if $a_1+(n-1)a_2\equiv 0\pmod2$, 
then $\Im(\f )
\subset \Ker (\sg)$, and $\f$ is not surjective.
Therefore, $a_1+(n-1)a_2\equiv 1\pmod2$.
If $a_1\equiv 1\pmod2$ and $a_2\equiv 0\pmod2$, then $\f = \nu_n^1$, and if
$a_1\equiv 1\pmod2$ and $a_2\equiv 1\pmod2$, then $n$ is odd and $\f = \nu_n^2$.
Assume $a_1\equiv0\pmod2$ and $a_2\equiv1\pmod2$. 
Then $n$ is even. Let $\alpha_1$ be the automorphism of $W$ determined by
$$\eqalign{
\alpha_1(s_1) &=\delta s_1,\cr
\alpha_1(s_i) &= s_i \quad {\rm for}\ i=2, \dots, n.\cr}
$$
Then
$$\eqalign{
(\alpha_1 \circ \f)(\sigma_1) = \alpha_1(t_1)&= c_1,\cr
(\alpha_1 \circ \f)(\sigma_i) = \alpha_1(t_i) &= c_{i-1}\cdot s_i \quad {\rm for}\ i=2, \dots, n, \cr}
$$
so $\alpha_1\circ\f=\nu_n^1$.
$ \square $

\nl
{\bf Lemma 3.5.} {\it For the Coxeter graph $B_6$, there is no extraordinary epimorphism $\f: A 
\to W$ such that $p \circ \f= \zeta_6$.

For the Coxeter graph $B_4$, there is no extraordinary epimorphism $\f: A \to W$ such that either 
$p \circ \f= \zeta_4^1$ or $p \circ \f= \zeta_4^2$.}

\nl
{\bf Proof.} The proof is a direct calculation. It can be carried out
in GAP or Magma by means of a Todd Coxeter enumeration leading to a
permutation representation of $W$, or by use of the package `Chevie' of
GAP. Here, we explain how we deal with $\zeta_6$. The other two cases
can be handled similarly.

We first compute the set $X_i= \{ xr_i; x \in C$ and $\order (xr_i) 
\neq 2\}$, for $i=2, \dots, 6$, where $r_i$ is as in Proposition 3.1.
We have $|X_i|=56$ for each $i$. For $i=3, \dots, 6$, let $Y_i$ be 
the set of $(i-1)$-tuples $(t_2, \dots, t_i)$ which satisfy

\smallskip
$\bullet$ $t_j \in X_j$, for $j=2, \dots, i$,

\smallskip
$\bullet$ $t_j t_{j+1} t_j = t_{j+1} t_j t_{j+1}$, for $j=2, \dots, i-1$,

\smallskip
$\bullet$ $t_j t_k = t_k t_j$, for $|j-k| \ge 2$.

\smallskip\noindent
We compute successively $Y_3$, $Y_4$, $Y_5$, and $Y_6$. We have
$|Y_3|=224$, $|Y_4|=192$, $|Y_5|=64$, and $Y_6=\emptyset$. The last
equality together with Lemma 2.1 proves the lemma. $\ \square$

\nl
{\bf Proof of Proposition 3.2.} Let $\f: A \to W$ be an extraordinary
epimorphism. By Proposition 3.1 and the fact that
$p(W) = \Sym_n$, we can assume that either $p \circ
\f$ is the standard epimorphism, or $n=6$ and $p \circ \f = \zeta_6$,
or $n=4$ and $p \circ \f = \zeta_4^i$, $i=1,2$. By Lemma 3.5, the last
two cases do not occur, so $p \circ \f$ is standard. Now,
%$$\displaylines{
%c_i= s_i s_{i-1} \cdots s_2 s_1 s_2 \cdots s_{i-1} s_i \quad{\rm for}\ i=2, \dots, n\cr
%(i-1,i) =s_i \quad {\rm for}\ i=2, \dots, n\cr
%[-1,1, \dots, 1] =s_1\cr
%[-1, -1, \dots, -1] = (s_1 s_2 \cdots s_n)^n\cr}
%$$
%thus
Proposition 3.2 follows from Lemma 3.4. $\ \square$

\bigskip
Let $A$ be an Artin group. By Lemma 2.3, in order to prove that $CA$
is a characteristic subgroup of $A$, it suffices to show that, for a
given automorphism $\Phi: A \to A$, the composition $\mu \circ \Phi: A
\to W$ with the standard epimorphism $\mu: A \to W$ is ordinary. For
example, for an Artin group of type $B_n$, we would like to prove that
$\mu \circ \Phi \neq \nu_n^i$, for $i=1,2$. Our general strategy is to
use the homomorphism $U: A \to \Perm( \Z \times T)$ defined in Lemma
3.6 below.

Let $T=\{ wsw^{-1}; w\in W\ {\rm and}\ s \in S\}$ be the set of
reflections in $W$. For all $i=1, \dots, n$, define the transformation
$u_i: \{ \pm 1 \} \times T \to \{ \pm 1 \} \times T$ by
$$
u_i(\e,t)= \left\{\matrix{ (\e, s_its_i)\hfill\quad&{\rm if}\ t \neq
s_i,\hfill\cr (-\e,t)\hfill\quad&{\rm if}\ t=s_i.\cr}
\right.$$
The following lemma can be found in [Bou, Ch.4, $n^o$1.4].

\nl
{\bf Lemma 3.6} (Bourbaki, [Bou]). {\it For each type $\Gamma$,
the assignment $s_i \mapsto u_i$ determines a homomorphism 
$u: W \to \Perm( \{ \pm 1 \} \times T)$, where $\Perm( \{ \pm 1 \} \times T)$ denotes the group of 
permutations of $ \{ \pm 1 \} \times T$. $\ \square$}

\bigskip
Now, for $i=1, \dots,n$, define $U_i: \Z \times T \to \Z \times T$ by
$$
U_i(k,t)=\left\{\matrix{
(k, s_its_i)\hfill\quad&{\rm if}\ t \neq s_i\hfill\cr
(k+1,t)\hfill\quad&{\rm if}\ t=s_i\cr}
\right.$$
Observe that each $U_i$ is a permutation of $\Z \times T$ with inverse
given by the same formula as for $U_i$, but with $k+1$ replaced by
$k-1$. The lemma below is easily derived by arguments similar to those
of [Bou, Ch.4, $n^o$1.4].

\nl
{\bf Lemma 3.7.} {\it For each type $\Gamma$,
the permutations $U_i$ of $\Z \times T$ satisfy the following properties.

\item{(1)} The assignment $\sigma_i \mapsto U_i$ 
determines a unique homomorphism $U: A \to
\Perm( \Z \times T)$.

\item{(2)} Let $g \in A$, and let $(k,r) \in \Z \times T$. Write $U(g)(k,r) =(k',r')$. Then
$$
u(\mu(g)) ((-1)^k,r) = ((-1)^{k'},r').
$$

\item{(3)} For $k\in \Z$ the element $L_k\in \Perm(\Z\times T)$ defined by
$L_k(l,t)= (k+l,t)$ commutes with each $U_i$. In particular,
for $g\in A$ we have $U(g)(k,t) = L_kU(g)(0,t)$.
$\quad\square$}

\bigskip
A direct consequence of the second part of this lemma is:

\nl
{\bf Corollary 3.8.} {\it Let $g_1,g_2 \in A$, and let $(k,r) \in \Z \times T$. Write $(k_1,r_1)= 
U(g_1)(k,r)$, and $(k_2,r_2)= U(g_2)(k,r)$. If $\mu(g_1)=\mu(g_2)$, then $r_1=r_2$ and $k_1 \equiv 
k_2 (\mod\ 2)$. $\quad\square$}

\bigskip
We now prove the second main result for $B_n$.

\nl
{\bf Proposition 3.9.} {\it Let $n \ge 3$. For the Coxeter graph $B_n$, the colored Artin group 
$CA$ is a characteristic subgroup of $A$.}

\bigskip
The key of the proof of Proposition 3.9 is the following lemma.

\nl
{\bf Lemma 3.10.} {\it Let $A$ be an Artin group of type $B_n$, $n\ge3$,
and let $\Phi: A \to A$ be a 
homomorphism. Then $\mu \circ \Phi \neq \nu_n^1$, and $\mu \circ \Phi \neq \nu_n^2$ (if $n$ is 
odd).}

\nl
{\bf Proof.} Assume that there exists a homomorphism $\Phi: A \to A$ such that $\mu \circ \Phi = 
\nu_n^i$, where $i\in \{ 1,2\}$. Consider the following elements of $T$.
\bigskip\centerline{\vbox{\halign{\hfill$#$&\ $#$\ &$#$\hfill \quad &\hfill$#$&\ $#$\ &$#$\hfill \quad 
&\hfill$#$&\ $#$\ &$#$\hfill\cr
r_1 &=& s_2,& r_2 &=& s_2s_3s_2, &r_3 &=& s_3,\cr
r_4 &=& s_1s_2s_1, &r_5 &=& s_1s_2s_3s_2s_1, &r_6 &=& s_2s_1s_2s_3s_2s_1s_2.\cr}}}
\nl
Put
$$
\tau_2= \sigma_1\sigma_2, \quad \hbox{ and }\quad \tau_3 =  \sigma_2 \sigma_1 \sigma_2 \sigma_3,
$$
so that
$\mu(\tau_j) = \nu_n^i(\sigma_j)$ for $j=2,3$.
Straightforward calculations give
\bigskip\centerline{\vbox{\halign{\hfill$#$&\ $#$\ &$#$\hfill \quad &\hfill$#$&\ $#$\ &$#$\hfill 
\quad &\hfill$#$&\ $#$\ &$#$\hfill\cr
U(\tau_2)(0,r_1) &=& (1,r_4), &U(\tau_2)(0,r_2) &=& (0,r_3), &U(\tau_2)(0,r_3) &=& (0,r_5),\cr
\noalign{\smallskip}
U(\tau_2)(0,r_4) &=& (0,r_1), &U(\tau_2)(0,r_5) &=& (0,r_6), &U(\tau_2)(0,r_6) &=& (0,r_2).\cr}}}
\nl
and
\bigskip\centerline{\vbox{\halign{\hfill$#$&\ $#$\ &$#$\hfill \quad &\hfill$#$&\ $#$\ &$#$\hfill 
\quad &\hfill$#$&\ $#$\ &$#$\hfill\cr
U(\tau_3)(0,r_1) &=& (0,r_2), &U(\tau_3)(0,r_2) &=& (1,r_4), &U(\tau_3)(0,r_3) &=& (1,r_6),\cr
\noalign{\smallskip}
U(\tau_3)(0,r_4) &=& (0,r_5), &U(\tau_3)(0,r_5) &=& (1,r_1), &U(\tau_3)(0,r_6) &=& (0,r_3).\cr}}}
\nl
Now $\mu(\tau_2) = \nu_n^i(\sigma_2) = \mu(\Phi(\sigma_2))$ and $\mu(\tau_3) = 
\nu_n^i(\sigma_3) = \mu (\Phi (\sigma_3))$, so, by Corollary 3.8, there exist integers
$a_1, \dots, a_6, b_1, \dots ,b_6$ such that
\bigskip\centerline{\vbox{\halign{\hfill$#$&\ $#$\ &$#$\hfill \quad &\hfill$#$&\ $#$\ &$#$\hfill 
\quad &\hfill$#$&\ $#$\ &$#$\hfill\cr
U(\Phi(\sigma_2))(0,r_1) &=& (2a_1+1,r_4),&U(\Phi(\sigma_2))(0,r_2) &=& (2a_2,r_3), 
&U(\Phi(\sigma_2))(0,r_3) &=& (2a_3,r_5),\cr
\noalign{\smallskip}
U(\Phi(\sigma_2))(0,r_4) &=& (2a_4,r_1), &U(\Phi(\sigma_2))(0,r_5) &=& (2a_5,r_6), 
&U(\Phi(\sigma_2))(0,r_6) &=& (2a_6,r_2).\cr}}}
\nl
and
\bigskip\centerline{\vbox{\halign{\hfill$#$&\ $#$\ &$#$\hfill \quad &\hfill$#$&\ $#$\ &$#$\hfill 
\quad &\hfill$#$&\ $#$\ &$#$\hfill\cr
U(\Phi(\sigma_3))(0,r_1) &=& (2b_1,r_2), &U(\Phi(\sigma_3))(0,r_2) &=& (2b_2+1,r_4), 
&U(\Phi(\sigma_3))(0,r_3) &=& (2b_3+1,r_6),\cr
\noalign{\smallskip}
U(\Phi(\sigma_3))(0,r_4) &=& (2b_4,r_5), &U(\Phi(\sigma_3))(0,r_5) &=& (2b_5+1,r_1), 
&U(\Phi(\sigma_3))(0,r_6) &=& (2b_6,r_3).\cr}}}
\nl
Using Lemma 3.7(3), we derive
\bigskip\centerline{\vbox{\halign{\hfill$#$&\ $#$\ &$#$\hfill \quad &\hfill$#$&\ $#$\ 
&$#$\hfill\cr 
U(\Phi(\sigma_2 \sigma_3 \sigma_2))(0,r_1) &=& (2a_1 +2a_5 +2b_4+1 ,r_6), &U(\Phi(\sigma_2 
\sigma_3 \sigma_2))(0,r_2) &=& (2a_2 +2a_6 +2b_3+1 ,r_2),\cr
\noalign{\smallskip}
U(\Phi(\sigma_2 \sigma_3 \sigma_2))(0,r_3) &=& (2a_1 +2a_3 +2b_5+2 ,r_4), &U(\Phi(\sigma_2 \sigma_3 
\sigma_2))(0,r_4) &=& (2a_2 +2a_4 +2b_1 ,r_3),\cr
\noalign{\smallskip}
U(\Phi(\sigma_2 \sigma_3 \sigma_2))(0,r_5) &=& (2a_3 +2a_5 +2b_6 ,r_5), &U(\Phi(\sigma_2 \sigma_3 
\sigma_2))(0,r_6) &=& (2a_4 +2a_6 +2b_2+1 ,r_1).\cr}}}
\nl
and
\bigskip\centerline{\vbox{\halign{\hfill$#$&\ $#$\ &$#$\hfill \quad &\hfill$#$&\ $#$\ 
&$#$\hfill\cr 
U(\Phi(\sigma_3 \sigma_2 \sigma_3)) (0,r_1) &=& (2a_2 +2b_1 +2b_3+1 ,r_6), &U(\Phi(\sigma_3 
\sigma_2 \sigma_3)) (0,r_2) &=& (2a_4+2b_1 +2b_2+1 ,r_2),\cr
\noalign{\smallskip}
U(\Phi(\sigma_3 \sigma_2 \sigma_3)) (0,r_3) &=& (2a_6 +2b_2 +2b_3+2 ,r_4), &U(\Phi(\sigma_3 
\sigma_2 \sigma_3)) (0,r_4) &=& (2a_5 +2b_4 +2b_6 ,r_3),\cr
\noalign{\smallskip}
U(\Phi(\sigma_3 \sigma_2 \sigma_3)) (0,r_5) &=& (2a_1 +2b_4 +2b_5+2 ,r_5), &U(\Phi(\sigma_3 
\sigma_2 \sigma_3)) (0,r_6) &=& (2a_3 +2b_5 +2b_6+1 ,r_1).\cr}}}
\nl
The equalities $U( \Phi (\sigma_2 \sigma_3 \sigma_2)) (0,r_i) = U( \Phi( \sigma_3 \sigma_2 
\sigma_3)) (0,r_i)$ for $i=2,3,4,5$ are equivalent to the equations
$$\eqalign{
a_2-a_4+a_6-b_1-b_2+b_3 &=0,\cr
a_1+a_3-a_6-b_2-b_3+b_5 &=0,\cr
a_2+a_4-a_5+b_1-b_4-b_6 &=0,\cr
-a_1+a_3+a_5-b_4-b_5+b_6 &=1.\cr}
$$
The sum of these four equations gives
$
2(a_2+a_3-b_2-b_4)=1
$
which clearly is impossible. $\ \square$

\nl
{\bf Proof of Proposition 3.9.} Let $\Delta= (\sigma_1 \cdots \sigma_n)^n \in A$, and recall that
$\delta= 
(s_1 \cdots s_n)^n \in W$. Let $\alpha_0: W \to W$ be the automorphism determined by
$$
\alpha_0(s_1)=s_1, \quad \alpha_0(s_i)= \delta s_i \quad {\rm for}\ i=2, \dots, n,
$$
and let $\Psi_0: A \to A$ be the homomorphism determined by
$$
\Psi_0(\sigma_1) =\sigma_1, \quad \Psi_0 (\sigma_i) = \Delta \sigma_i \quad {\rm for}\ i=2, \dots, n.
$$
Then $\mu \circ \Psi_0= \alpha_0 \circ \mu$. For $n$ even, let $\alpha_1: W \to W$ be the 
automorphism determined by
$$
\alpha_1(s_1) = \delta s_1, \quad \alpha_1 (s_i) = s_i \quad {\rm for}\ i=2, \dots, n,
$$
and let $\Psi_1: A \to A$ be the homomorphism determined by
$$
\Psi_1 (\sigma_1) = \Delta \sigma_1, \quad \Psi_1 (\sigma_i) = \sigma_i \quad {\rm for}\ i=2, \dots, n.
$$
Then $\mu \circ \Psi_1 = \alpha_1 \circ \mu$. 
Observe also that, if $w \in W$, then $\Conj_w \circ \mu= \mu \circ \Conj_g$, where $g\in A$ is such
that $\mu(g)=w$.
By [Fra], $\Aut (W)$ is generated as a monoid by 
$\{\Conj_w; w \in W\} \cup \{ \alpha_0, \alpha_1 \}$ if $n$ is even, and by $\{\Conj_w; w \in W\} 
\cup \{ \alpha_0\}$ if $n$ is odd. In particular, by the above observations,  
if $\alpha$ is an automorphism of $W$, then there exists a homomorphism $\Psi: A \to 
A$ such that $\mu \circ \Psi = \alpha \circ \mu$.

Now, let $\Phi: A \to A$ be an automorphism. In order to prove
Proposition 3.9, it suffices to show that $\mu \circ \Phi$ cannot be
extraordinary.  Suppose that $\mu \circ \Phi$ is extraordinary. By
Proposition 3.2, there exist an automorphism $\alpha: W \to W$ and $i
\in \{1,2\}$ such that $\mu \circ \Phi = \alpha \circ \nu_n^i$. By the
preceding argument, one can find a homomorphism $\Psi: A \to A$ such
that $\mu \circ \Psi =
\alpha^{-1} \circ \mu$, so $\mu \circ (\Psi \circ \Phi)= \nu_n^i$. 
%(Note that
%the function $\Phi$ of Lemma 3.10 is not 
%necessarily assumed to be an automorphism.) 
But this equality contradicts Lemma 3.10. $\ \square$

%%%%%%%%%%%%%%%%%%%%%%%%%%%%%%%%%%%%%%%%%%%%%%%%%%%%%%%%%
\bigskip\nl
{\titre 4. Artin groups of type $D_n$}

\bigskip
The strategy for classifying the epimorphisms $A \to W$ for the Artin groups of type $D_n$, $n 
\ge 4$, is the same as the strategy for the Artin groups of type $B_n$. First, we classify the 
epimorphisms $A \to \Sym_n$ up to conjugation. This has also been
carried out by Zinde [Zin], and the answer is as follows.
Label the nodes of the Coxeter graph $D_n$ as follows:
$$\Dn$$
Define the standard epimorphism $\eta: A \to \Sym_n$ as the homomorphism determined by
$$\eqalign{
&\eta(\sigma_1)  = \eta(\sigma_2) = (1,2),\cr
&\eta(\sigma_i) = (i-1,i) \quad {\rm for}\ i=3, \dots, n.\cr}
$$

\noindent
{\bf Proposition 4.1} (Zinde, [Zin]). {\it Let $n \ge 4$. For the Coxeter graph $D_n$, every 
epimorphism $A \to \Sym_n$ is conjugate to $\eta$ if $n\neq 4$ and $n \neq 6$.

If $n=6$, then there are precisely two conjugacy classes of epimorphisms. These are
represented by $\eta, \zeta_6: A 
\to \Sym_6$, respectively, where
\bigskip\centerline{\vbox{\halign{\hfill$#$&\ $#$\ &$#$\hfill\quad &\hfill$#$&\ $#$\ &$#$\hfill 
\quad &\hfill$#$&\ $#$\ &$#$\hfill\cr
\zeta_6(\sigma_1) &=& (1,2)(3,4)(5,6), &\zeta_6(\sigma_2) &=& (1,2)(3,4)(5,6), &\zeta_6(\sigma_3) 
&=& (2,3)(1,5)(4,6),\cr
\noalign{\smallskip}
\zeta_6(\sigma_4) &=& (3,1)(2,4)(5,6), &\zeta_6(\sigma_5) &=& (1,2)(3,5)(4,6), &\zeta_6(\sigma_6) 
&=& (2,3)(1,4)(5,6).\cr}}}
\bigskip

If $n=4$, for each epimorphism $\f: A \to \Sym_4$, there is a diagram automorphism $\Upsilon$ of 
$A$ such that $\f \circ \Upsilon$ is conjugate to one of the three epimorphisms $\eta, \zeta_4^1, 
\zeta_4^2: A\to \Sym_4$, where
$$
\zeta_4^1(\sigma_1) = \zeta_4^1(\sigma_2) = \zeta_4^1(\sigma_4) = (1,2,3,4), \quad 
\zeta_4^1(\sigma_3) = (2,1,3,4),
$$
and}
$$
\zeta_4^2(\sigma_1) = \zeta_4^2(\sigma_2) = (1,2,3,4), \quad \zeta_4^2(\sigma_3)= (2,1,3,4), 
\quad \zeta_4^2(\sigma_4) = (4,3,2,1). \quad \square
$$
 
It is well known that $W$ embeds into the Weyl group of type $B_n$
with image an index two subgroup.  The embedding can be made explicit
by means of the isomorphism of the Weyl group of type $B_n$ with
$(\Z/2\Z)^n\semidirect\Sym_n$ described in the discussion preceding
Proposition 3.2. The embedding maps $s_1$ to $(1,1,0,\ldots,0) \cdot
(1,2)$ and maps $s_i$ to $(i-1,i)$ for $i\ge2$.  So, if we denote by
$H$ the index 2 subgroup of $(\Z/2\Z)^n$ of all elements
$(a_1,\ldots,a_n)$ with $a_1+\cdots+a_n=0$, we find that the image is
$H\semidirect \Sym_n$.  Let $\delta= (s_1s_2) (s_3s_1s_2s_3) \cdots
(s_n \cdots s_3s_1s_2s_3\cdots s_n)$ denote the longest element of
$W$. Then $\delta $ is mappped to $(1,1,\ldots,1)\in H$ if $n$ is even and
to $(0,1,\ldots,1)\in H$ if $n$ is odd. For $n$ odd, write
$u_1=\delta$ and $u_i = s_is_{i-1}\cdots s_2u_1s_2\cdots s_{i-1}s_i$
if $i>1$.  We use Proposition 4.1 to prove Theorem 1.3 for $D_n$.

\nl
{\bf Proposition 4.2.} {\it 
Let $n \ge 4$. For
the Coxeter graph $D_n$,
there are 
no
extraordinary epimorphisms $A \to W$ if $n$ is even.
If $n$ is odd, there are precisely two
equivalence classes of 
extraordinary epimorphisms $A \to W$. These are
represented by $\nu_n^1, \nu_n^2: A \to W$, where
$$\displaylines{
\nu_n^1(\sigma_1) = \nu_n^1(\sigma_2) = s_2 \delta = u_2s_2, \cr
\nu_n^1(\sigma_i) = %s_i s_{i-1} \cdots s_2 \delta s_2 \cdots s_{i-1} = 
u_is_i \quad {\it for}\ i=2, \dots, 
n.\cr}
$$}
{\it and}
$$\displaylines{
\nu_n^2(\sigma_1)= u_1 s_2,\cr
\nu_n^2(\sigma_i) = % s_i s_{i-1} \cdots s_2 \delta s_2 \cdots s_{i-1}
u_is_i \quad {\it for}\ i=2, \dots, 
n.\cr}
$$

The following two lemmas are preparations to the proof of Proposition 4.2.
Denote by $p: W \to \Sym_n$ the natural projection with kernel $H$.

\nl
{\bf Lemma 4.3.} {\it Let $\f: A \to W$ be an extraordinary
epimorphism such that $p \circ \f$ is the standard epimorphism $\eta:
A \to \Sym_n$. Then $n$ is odd and $\f$ is equivalent to $\nu_n^1$ or $\nu_n^2$.}

\nl
{\bf Proof.} Let $t_i= \f(\sigma_i)$ for $i=1,2, \dots, n$. Write
$t_1= v_1 \cdot (1,2)$, and $t_i=v_i \cdot (i-1,i)$ (for $i=2, \dots,
n$), where $v_i \in H$. The same argument as in the proof of Lemma 3.3
shows that there exist $a, a_1, a_2, \dots, a_n \in \Z/2\Z$ such that
$v_1=(a_1,1+a_1, a, \dots, a)$, and $v_i=(a, \dots, a, a_i,
1+a_i, a, \dots, a)$ ($a_i$ being in the $(i-1)$-st entry) for
$i=2, \dots, n$. The condition $v_i \in H$ implies that $n$ is odd and
$a=1$. Using again the same argument as in the proof of Lemma 3.3,
we see that $a_i=1$ for all $i=2, \dots, n$ (not
necessarily for $i=1$), up to an automorphism of $W$. $\ \square$

\nl
{\bf Lemma 4.4.} {\it For the Coxeter graph $D_6$, there is no extraordinary epimorphism $\f: A 
\to W$ such that $p \circ \f= \zeta_6$.

For the Coxeter graph $D_4$, there is no extraordinary epimorphism $\f: A \to W$ such that either 
$p \circ \f = \zeta_4^1$ or $p \circ \f= \zeta_4^2$.}

\nl
{\bf Proof.} The proof is a straightforward calculation
which can be carried out in 
GAP (for instance with the package `Chevie' of GAP) or Magma. It follows the 
same scheme as the proof of Lemma 3.5. $\ \square$

\nl
{\bf Proof of Proposition 4.2.} Let $\f: A \to W$ be an extraordinary epimorphism. By Proposition 
4.1, we can assume that either $p \circ \f$ is the standard epimorphism, or $n=6$ and $p \circ 
\f= \zeta_6$, or $n=4$ and $p \circ \f \circ \Upsilon =\zeta_4^i$, where $\Upsilon$ is a diagram 
automorphism of $A$, and $i\in \{1,2\}$. By Lemma 4.4, the last two cases cannot happen, so $p 
\circ \f$ is standard.
Proposition 4.2 now follows from Lemma 4.3. $\ \square$

\bigskip
We turn now to the proof of the second main result for $D_n$.

\nl
{\bf Proposition 4.5.} {\it Let $n \ge 4$. For the Coxeter graph
$D_n$, the colored Artin group $CA$ is a characteristic subgroup of
$A$.}

\nl
{\bf Remark.} By Artin's Theorem (1.1), for the Coxeter graph $A_n$,
each epimorphism $A \to W$ is standard, except for $n=3$. In
particular, this shows that, for an Artin group of type $A_n$, $n \neq
3$, the colored Artin group $CA$ is a characteristic subgroup of
$A$. Now, the proof given below also applies to $A_3$, viewed as the
Coxeter graph $D_3$, and so provides an alternative proof to Theorem
1.2 for $n=4$.

\bigskip
The following lemma is a preparation to the proof of Proposition 4.5.

\nl
{\bf Lemma 4.6.} {\it Assume $n$ is odd. Let $A$ be an Artin group of type $D_n$, and let $\Phi: 
A \to A$ be a homomorphism. Then $\mu \circ \Phi \neq \nu_n^1$, and $\mu \circ \Phi \neq 
\nu_n^2$.}

\nl
{\bf Proof.} Assume that there exists a homomorphism $\Phi: A \to A$ such that $\mu \circ \Phi= 
\nu_n^i$, where $i\in \{ 1,2\}$. Consider the following six reflections.
\bigskip\centerline{\vbox{\halign{\hfill$#$&\ $#$\ &$#$\hfill\quad &\hfill$#$&\ $#$\ &$#$\hfill 
\quad &\hfill$#$&\ $#$\ &$#$\hfill\cr
r_1 &=& s_1, &r_2 &=& s_3s_1s_3, &r_3 &=& s_2s_3s_1s_3s_2,\cr
r_4 &=& s_3, &r_5 &=& s_2s_3s_2, &r_6 &=& s_2.\cr}}}
\nl
Let
$$
\Delta= (\sigma_1 \sigma_2) (\sigma_3 \sigma_1 \sigma_2 \sigma_3) \cdots (\sigma_n \cdots \sigma_3 
\sigma_1 \sigma_2 \sigma_3 \cdots \sigma_n).
$$
Recall the map $U$ defined in Lemma 3.7.  The fact that
$\delta=(s_1s_2) (s_3s_1s_2s_3) \cdots (s_n \cdots s_3s_1s_2s_3 \cdots
s_n)$ is a reduced expression for the longest element of $W$ implies
that $U(\Delta)(0,r_k)= (1,\delta r_k\delta)$ for $k=1, \dots,
6$. Furthermore, since $n$ is odd, $\delta s_1 \delta = s_2$, $\delta
s_2 \delta = s_1$, and $\delta s_j \delta = s_j$ for $j=3, \dots,
n$. Let
$$
\tau_2 = \sigma_2 \Delta, \quad \tau_3=\sigma_3 \sigma_2 \Delta \sigma_2,
$$
so that
$\mu(\tau_2) = s_2 \delta =\nu_n^i(\sigma_2)$ and $\mu(\tau_3) = 
s_3s_2 \delta s_2 = \nu_n^i(\sigma_3)$.
A straightforward calculation gives
\bigskip\centerline{\vbox{\halign{\hfill$#$&\ $#$\ &$#$\hfill\quad &\hfill$#$&\ $#$\ &$#$\hfill 
\quad &\hfill$#$&\ $#$\ &$#$\hfill\cr
U(\tau_2)(0,r_1) &=& (2,r_6), &U(\tau_2)(0,r_2) &=& (1,r_4), &U(\tau_2)(0,r_3) &=& (1,r_2),\cr
\noalign{\smallskip}
U(\tau_2)(0,r_4) &=& (1,r_5), &U(\tau_2)(0,r_5) &=& (1,r_3), &U(\tau_2)(0,r_6) &=& (1,r_1).\cr}}}
\nl
and
\bigskip\centerline{\vbox{\halign{\hfill$#$&\ $#$\ &$#$\hfill\quad &\hfill$#$&\ $#$\ &$#$\hfill 
\quad &\hfill$#$&\ $#$\ &$#$\hfill\cr
U(\tau_3)(0,r_1) &=& (2,r_5), &U(\tau_3)(0,r_2) &=& (1,r_1), &U(\tau_3)(0,r_3) &=& (2,r_4),\cr
\noalign{\smallskip}
U(\tau_3)(0,r_4) &=& (1,r_3), &U(\tau_3)(0,r_5) &=& (1,r_6), &U(\tau_3)(0,r_6) &=& (2,r_2).\cr}}}
\nl
Now $\mu(\tau_2) = \nu_n^i(\sigma_2) = \mu(\Phi(\sigma_2))$ and $\mu(\tau_3) = 
\nu_n^i(\sigma_3)= \mu(\Phi(\sigma_3))$ so, by Corollary 3.8, there exist 
integers $a_1, \dots, a_6, b_1, \dots, b_6$ such that
\bigskip\centerline{\vbox{\halign{\hfill$#$&\ $#$\ &$#$\hfill\quad &\hfill$#$&\ $#$\ &$#$\hfill 
\quad &\hfill$#$&\ $#$\ &$#$\hfill\cr
U(\Phi(\sigma_2))(0,r_1) &=& (2a_1,r_6), &U(\Phi(\sigma_2))(0,r_2) &=& (2a_2+1,r_4), 
&U(\Phi(\sigma_2))(0,r_3) &=& (2a_3+1,r_2),\cr
\noalign{\smallskip}
U(\Phi(\sigma_2))(0,r_4) &=& (2a_4+1,r_5), &U(\Phi(\sigma_2))(0,r_5) &=& (2a_5+1,r_3), 
&U(\Phi(\sigma_2))(0,r_6) &=& (2a_6+1,r_1).\cr}}}
\nl
and
\bigskip\centerline{\vbox{\halign{\hfill$#$&\ $#$\ &$#$\hfill\quad &\hfill$#$&\ $#$\ &$#$\hfill 
\quad &\hfill$#$&\ $#$\ &$#$\hfill\cr
U(\Phi(\sigma_3))(0,r_1) &=& (2b_1,r_5), &U(\Phi(\sigma_3))(0,r_2) &=& (2b_2+1,r_1),  
&U(\Phi(\sigma_3))(0,r_3) &=& (2b_3,r_4),\cr
\noalign{\smallskip}
U(\Phi(\sigma_3))(0,r_4) &=& (2b_4+1,r_3), &U(\Phi(\sigma_3))(0,r_5) &=& (2b_5+1,r_6), 
&U(\Phi(\sigma_3))(0,r_6) &=& (2b_6,r_2).\cr}}}
\nl
We find that
\bigskip\centerline{\vbox{\halign{\hfill$#$&\ $#$\ &$#$\hfill\quad &\hfill$#$&\ $#$\ 
&$#$\hfill\cr
U(\Phi(\sigma_2\sigma_3\sigma_2))(0,r_1) &=& (2a_1+2a_2+2b_6+1,r_4), 
&U(\Phi(\sigma_2\sigma_3\sigma_2))(0,r_2) &=& (2a_2+2a_3+2b_4+3,r_2),\cr
\noalign{\smallskip}
U(\Phi(\sigma_2\sigma_3\sigma_2))(0,r_3) &=& (2a_1+2a_3+2b_2+2,r_6), 
&U(\Phi(\sigma_2\sigma_3\sigma_2))(0,r_4) &=& (2a_4+2a_6+2b_5+3,r_1),\cr
\noalign{\smallskip}
U(\Phi(\sigma_2\sigma_3\sigma_2))(0,r_5) &=& (2a_4+2a_5+2b_3+2,r_5), 
&U(\Phi(\sigma_2\sigma_3\sigma_2))(0,r_6) &=& (2a_5+2a_6+2b_1+2,r_3).\cr}}}
\nl
and
\bigskip\centerline{\vbox{\halign{\hfill$#$&\ $#$\ &$#$\hfill\quad &\hfill$#$&\ $#$\ 
&$#$\hfill\cr
U(\Phi(\sigma_3\sigma_2\sigma_3))(0,r_1) &=& (2a_5+2b_1+2b_3+1,r_4), 
&U(\Phi(\sigma_3\sigma_2\sigma_3))(0,r_2) &=& (2a_1+2b_2+2b_6+1,r_2),\cr
\noalign{\smallskip}
U(\Phi(\sigma_3\sigma_2\sigma_3))(0,r_3) &=& (2a_4+2b_3+2b_5+2,r_6), 
&U(\Phi(\sigma_3\sigma_2\sigma_3))(0,r_4) &=& (2a_3+2b_2+2b_4+3,r_1),\cr
\noalign{\smallskip}
U(\Phi(\sigma_3\sigma_2\sigma_3))(0,r_5) &=& (2a_6+2b_1+2b_5+2,r_5), 
&U(\Phi(\sigma_3\sigma_2\sigma_3))(0,r_6) &=& (2a_2+2b_4+2b_6+2,r_3).\cr}}}
\nl
The equalities $U(\Phi(\sigma_2 \sigma_3 \sigma_2))(0, r_i)= U(\Phi( \sigma_3 \sigma_2 
\sigma_3))(0,r_i)$ for $i=1,2,4,5$ are equivalent to the equations
$$\eqalign{
a_1+a_2-a_5-b_1-b_3+b_6 &=0\cr
-a_1+a_2+a_3-b_2+b_4-b_6 &=-1\cr
-a_3+a_4+a_6-b_2-b_4+b_5 &=0\cr
a_4+a_5-a_6-b_1+b_3-b_5 &=0\cr}
$$
The sum of these four equations gives
$2(a_2+a_4-b_1-b_2)=-1$, a contradiction. $\ \square$

\nl
{\bf Proof of Proposition 4.5.} Let $\Phi: A \to A$ be an automorphism.
Suppose that $\mu \circ \Phi$ is extraordinary. By Proposition 4.2, $n$ is odd, and there exist 
an automorphism $\alpha: W \to W$ and $i \in \{1,2\}$ such that $\mu \circ \Phi= \alpha \circ 
\nu_n^i$. By [Fra], $\Aut W =\{\Conj_w; w\in W\}$. Let $w \in W$ be such that $\alpha^{-
1}=\Conj_w$, and take $g \in A$ such that $\mu(g)=w$. Then $\alpha^{-1} \circ \mu = \Conj_w \circ 
\mu=\mu \circ \Conj_g$, whence $\mu \circ (\Conj_g \circ \Phi)= \nu_n^i$. This equality 
contradicts Lemma 4.6.  Hence $\mu\circ\Phi$ is ordinary, and so the
proposition follows from Lemma 2.3.  $\ \square$

%%%%%%%%%%%%%%%%%%%%%%%%%%%%%%%%%%%%%%%%%%%%%%%%%%%%%%%%%%
\bigskip\nl
{\titre 5. The remainder}

\bigskip
In this section, we classify up to equivalence the epimorphisms $A \to
W$, and we prove that $CA$ is a characteristic subgroup of $A$, for
the Artin groups of type $H_3$, $H_4$, $F_4$, $E_6$, $E_7$,
$E_8$. This finishes our proofs of Theorems 1.3 and 1.4.  We start
with the Artin group of type $H_3$, which requires a separate
treatment, because it is the only Artin group among those of type
$H_3$, $H_4$, $F_4$, $E_6$, $E_7$, $E_8$ which has extraordinary
epimorphisms.  For the corresponding Coxeter diagrams, we adopt the
labelings
$$
{\vtriple{\scriptstyle1}\over\circ\over{}
\kern-4pt{\buildrel {\scriptstyle 5}\over\lijntje}\kern-4pt\vtriple{\scriptstyle{2}}\over\circ\over{}
\kern-1pt\lijntje\kern-1pt
\vtriple{\scriptstyle 3}\over\circ\over{}
}
\qquad\hbox{ and }\qquad
{\vtriple{\scriptstyle1}\over\circ\over{}
\kern-4pt{\buildrel {\scriptstyle 5}\over\lijntje}\kern-4pt\vtriple{\scriptstyle{2}}\over\circ\over{}
\kern-1pt\lijntje\kern-1pt
\vtriple{\scriptstyle 3}\over\circ\over{}
\kern-1pt\lijntje\kern-1pt
\vtriple{\scriptstyle 4}\over\circ\over{}\kern-1pt}
$$
for $H_3$ and $H_4$ and the labeling of [Bou] for the remaining types.

\nl
{\bf Lemma 5.1.} {\it For the Coxeter graph $H_3$, there are precisely
two conjugacy classes of 
extraordinary epimorphisms $A \to W$.
These are represented by $\nu_3^1, \nu_3^2: A \to W$, where
$$
\nu_3^1(\sigma_1) = s_1s_2s_3, \quad \nu_3^1(\sigma_2) = s_2s_1s_2s_3s_2, \quad \nu_3^1(\sigma_3) 
= s_3s_2s_1,
$$
and}
$$\displaylines{
\nu_3^2(\sigma_1) = s_1s_2s_1 s_2s_3s_2 s_1s_2s_3, \quad \nu_3^2(\sigma_2) = s_2s_1s_2 s_3s_2s_1 
s_2s_1s_3 s_2s_1s_2 s_3,\quad \nu_3^2(\sigma_3) = s_2s_3s_2 s_1s_2s_1 s_3s_2s_1.}
$$

\noindent
{\bf Proof.} The proof is a straightforward calculation, which can be
made in GAP (for instance, with the package `Chevie') or Magma. First,
we calculate a representative for each conjugacy class of $W$. Let
$C_1$ be the set of these representatives. We compute $C_2= \{ x \in
C_1; \sg(x)=-1$ and $\order (x) \neq 2 \}$.  We have $|C_2|=3$. Let
$$\eqalign{
X_1 &= \{ (t_1,t_2) \in C_2 \times W\ ;\ t_1 \neq t_2,\ t_1t_2t_1t_2t_1 = t_2t_1t_2t_1t_2 \}, \cr
X_2 &= \{ (t_1,t_2,t_3) \in X_1 \times W\ ;\ t_1t_3 = t_3t_1,\ t_2t_3t_2 = t_3t_2t_3 \}. \cr}
$$
We successively compute $X_1$ and $X_2$. We have $|X_1|=16$, and $|X_2|=10$. We verify that, if 
$(t_1,t_2,t_3) \in X_2$, then $\{ t_1,t_2,t_3 \}$ generates $W$. Let $\sim$ be the equivalence 
relation on $X_2$ defined by $(t_1,t_2,t_3) \sim (t_1',t_2',t_3')$ if there exists $w \in W$ such 
that $t_i'=wt_iw^{-1}$ for all $i=1,2,3$. We compute a representative for each equivalence class. 
Let $Y$ be the set of these representatives. Then
$$\displaylines{
\qquad Y=\{ (s_1s_2s_3, s_2s_1s_2s_3s_2, s_3s_2s_1), \hfill\cr
\hfill (s_1s_2s_1 s_2s_3s_2 s_1s_2s_3,  s_2s_1s_2 s_3s_2s_1 s_2s_1s_3 s_2s_1s_2 s_3, s_2s_3s_2 
s_1s_2s_1 s_3s_2s_1) \}.\quad \square \qquad\cr}
$$

\noindent
{\bf Proposition 5.2.} {\it For the Coxeter graph $H_3$, there is a unique equivalence class of 
extraordinary epimorphisms $A \to W$. It is represented by $\nu_3^1: A \to W$.}

\nl
{\bf Proof.} Let $\alpha: W \to W$ be the automorphism determined by
$$
\alpha(s_1) = s_1s_2s_3 s_2s_1, \quad \alpha(s_2) = s_2s_1s_2s_3s_2s_1s_2, \quad \alpha(s_3) = 
s_1s_2s_1 s_2s_1.
$$
Then $\alpha \circ \nu_3^1 = \nu_3^2$. $\ \square$

\nl
{\bf Proposition 5.3.} {\it For the Coxeter graph $H_3$, the colored Artin group $CA$ is a 
characteristic subgroup of $A$.}

\nl
{\bf Proof.} Let $\Phi: A \to A$ be an automorphism. In order to prove
Proposition 5.3, it suffices to show that $\mu \circ \Phi$
is an ordinary epimorphism, cf.\ Lemma 2.3. In particular, we need
only establish the claim that $\mu \circ \Phi \neq
\nu_3^1$ and $\mu \circ \Phi \neq \nu_3^2$.
For, if $\mu\circ\Phi$ is extraordinary, then by Lemma 5.1 there are
$w\in W$ and $i\in\{1,2\}$ such that $\mu\circ\Phi =
\Conj_w\circ\nu_3^i$; taking $g\in A$ with $\mu(g^{-1}) = w$, we find
$\mu\circ(\Conj_{g}\circ\Phi) = \Conj_{w^{-1}}\circ\mu\circ\Phi = \nu_3^i$,
so we can apply the claim with
$\Conj_{g}\circ\Phi$ instead of $\Phi$ to derive a contradiction.

In order to settle the claim,
suppose that $\mu \circ \Phi= \nu_3^1$. Let
\bigskip\centerline{\vbox{\halign{\hfill$#$&\ $#$\ &$#$\hfill\quad&\hfill$#$&\ $#$\ 
&$#$\hfill\quad&\hfill$#$&\ $#$\ &$#$\hfill\cr
r_1 &=& s_1, &r_2 &=& s_2, &r_3 &=& s_3,\cr
r_4 &=& s_1s_2s_1, &r_5 &=& s_2s_3s_2, &r_6 &=& s_2s_1s_2,\cr
r_7 &=& s_1s_2s_3 s_2s_1, &r_8 &=& s_3s_2s_1s_2s_3, &r_9 &=& s_1s_2s_1 s_2s_1,\cr
r_{10} &=& s_1s_3s_2 s_1s_2s_3 s_1, &r_{11} &=& s_2s_1s_2 s_3s_2s_1 s_2, &r_{12} &=& s_2s_1s_3 
s_2s_1s_2 s_3s_1s_2,\cr
r_{13} &=& s_1s_2s_1 s_2s_3s_2 s_1s_2s_1, &r_{14} &=& s_1s_2s_1 s_3s_2s_1 s_2s_3s_1 s_2s_1, 
&r_{15} &=& s_2s_1s_2 s_1s_3s_2 s_1s_2s_3 s_1s_2s_1 s_2.\cr}}}
\nl
It is easily verified that $T= \{wsw^{-1}; w\in W, s\in S\}= \{r_1, \dots, r_{15} \}$. Let
$$
\tau_2= \sigma_2 \sigma_1 \sigma_2 \sigma_3 \sigma_2, \quad \tau_3= \sigma_3 \sigma_2 \sigma_1.
$$
A direct calculation gives
\bigskip\centerline{\vbox{\halign{\hfill$#$&\ $#$\ &$#$\hfill\quad&\hfill$#$&\ $#$\ 
&$#$\hfill\quad&\hfill$#$&\ $#$\ &$#$\hfill\cr
U(\tau_2)(0,r_1) &=& (0,r_{12}), &U(\tau_2)(0,r_2) &=& (1,r_5), &U(\tau_2)(0,r_3) &=& (1,r_9),\cr
U(\tau_2)(0,r_4) &=& (0,r_{15}), &U(\tau_2)(0,r_5) &=& (1,r_{11}), &U(\tau_2)(0,r_6) &=& 
(0,r_4),\cr
U(\tau_2)(0,r_7) &=& (0,r_8), &U(\tau_2)(0,r_8) &=& (1,r_6), &U(\tau_2)(0,r_9) &=& (0,r_{13}),\cr
U(\tau_2)(0,r_{10}) &=& (0,r_3), &U(\tau_2)(0,r_{11}) &=& (0,r_1), &U(\tau_2)(0,r_{12}) &=& 
(1,r_2),\cr
U(\tau_2)(0,r_{13}) &=& (0,r_{14}), &U(\tau_2)(0,r_{14}) &=& (0,r_{10}), &U(\tau_2)(0,r_{15}) &=& 
(0,r_7).\cr}}}
\nl
and
\bigskip\centerline{\vbox{\halign{\hfill$#$&\ $#$\ &$#$\hfill\quad&\hfill$#$&\ $#$\ 
&$#$\hfill\quad&\hfill$#$&\ $#$\ &$#$\hfill\cr
U(\tau_3)(0,r_1) &=& (1,r_8), &U(\tau_3)(0,r_2) &=& (0,r_{10}), &U(\tau_3)(0,r_3) &=& (0,r_2),\cr
U(\tau_3)(0,r_4) &=& (1,r_5), &U(\tau_3)(0,r_5) &=& (0,r_{12}), &U(\tau_3)(0,r_6) &=& (0,r_7),\cr
U(\tau_3)(0,r_7) &=& (1,r_3), &U(\tau_3)(0,r_8) &=& (0,r_{11}), &U(\tau_3)(0,r_9) &=& (0,r_1),\cr
U(\tau_3)(0,r_{10}) &=& (0,r_6), &U(\tau_3)(0,r_{11}) &=& (0,r_{14}), &U(\tau_3)(0,r_{12}) &=& 
(0,r_{15}),\cr
U(\tau_3)(0,r_{13}) &=& (0,r_4), &U(\tau_3)(0,r_{14}) &=& (0,r_9), &U(\tau_3)(0,r_{15}) &=& 
(0,r_{13}).\cr}}}
\nl
Since $\mu( \Phi( \sigma_2)) = \nu_3^1(\sigma_2)= s_2s_1s_2 s_3s_2 =\mu(\tau_2)$, and $\mu( \Phi( 
\sigma_3))= \nu_3^2(\sigma_3)= s_3s_2s_1 =\mu(\tau_3)$, by Corollary 3.8, there exist $b_1, 
\dots, b_{15}, c_1, \dots, c_{15} \in \Z$ such that

\smallskip\noindent
$\bullet$ if $U(\tau_2) (0,r_i)= (k_i,r_j)$, then $U(\Phi( \sigma_2)) (0,r_i)= (2b_i+k_i,r_j)$;

\smallskip\noindent
$\bullet$ if $U(\tau_3) (0,r_i)= (k_i',r_j')$, then $U(\Phi( \sigma_3)) (0,r_i)= 
(2c_i+k_i',r_j')$.

\smallskip\noindent
Now, the equalities $U(\Phi( \sigma_2 \sigma_3 \sigma_2)) (0,r_i)= U(\Phi( \sigma_3 \sigma_2 
\sigma_3)) (0,r_i)$ for $i=2,5,6, 8,10,11$, are equivalent to the equations
$$\eqalign{
b_2-b_{10}+b_{12}-c_2-c_3+c_5 &=-1,\cr
b_5-b_{12}+b_{14}-c_2-c_5+c_{11} &=0,\cr
b_5+b_6-b_7+c_4-c_6-c_8 &=-1,\cr
b_7+b_8-b_{11}-c_1+c_6-c_8 &=0,\cr
b_2-b_6+b_{10}+c_3-c_4-c_{10} &=0,\cr
b_8+b_{11}-b_{14}+c_1-c_{10}-c_{11} &=-1.\cr}
$$
The sum of these six equations gives
$2(b_2+b_5+b_8-c_2-c_8-c_{10})=-3$
which clearly cannot hold. So, $\mu \circ \Phi \neq \nu_3^1$.

Suppose $\mu \circ \Phi= \nu_3^2$. Using the same procedure as before, we find that the 
equalities $U(\Phi (\sigma_2 \sigma_3 \sigma_2)) (0,r_i)$
$= U(\Phi( \sigma_3 \sigma_2 
\sigma_3))(0,r_i)$ for $i=2,4,7, 8,10,12$, are equivalent to the equations
$$\eqalign{
b_2-b_7+b_{15}-c_2+c_4-c_{13} &=0,\cr
b_2+b_4-b_{15}-c_4+c_6-c_{12} &=0,\cr
b_7-b_{10}+b_{12}-c_7+c_{13}-c_{14} &=0,\cr
b_3+b_8-b_9-c_2-c_8+c_{10} &=-1,\cr
-b_3+b_8+b_{10}-c_7-c_{10}+c_{14} &=0,\cr
-b_4+b_9+b_{12}-c_6+c_8-c_{12} &=0.\cr}
$$
The sum of these six equations gives
$2(b_2+b_8+b_{12}-c_2-c_7-c_{12}) = -1$,
a contradiction. So, $\mu \circ \Phi \neq \nu_3^2$. $\ \square$

\bigskip
Now, we study the remaining cases.

\nl
{\bf Proposition 5.4.} {\it For each of the Coxeter graphs $H_4$, $F_4$, $E_6$, $E_7$, $E_8$, all 
epimorphisms $A \to W$ are ordinary.}

\nl
{\bf Proof.}  Suppose $\f :A\to W$ is a morphism of groups.  In a
case-by-case analysis, we find all possible morphisms by listing
all possible $
(\f(\sigma_1),\ldots,\f(\sigma_n))$ up to conjugacy in $W$.  This
is done by computing all sequences $(x_1,\ldots,x_n)\in W^n$
satisfying $ [x_i, x_j
\rangle^{m_{i,j}} = [x_j, x_i \rangle^{m_{i,j}} $ for all $i,j\in
\{1,\ldots,n\}$, up to conjugacy in
$W$, and next verifying whether the $x_i$ generate $W$.  Excluding
the ordinary and trivial cases, we can impose some a priori conditions
on the $x_i$. For instance, in view of Lemma 2.1, we may assume that
some $x_i$ has order at least 3, and surjectivity of $\f$ implies
that some $x_i$ has sign $-1$.

It is well known (and easy to verify) that $x_i$ and $x_j$ are
conjugate whenever $i$ and $j$ are joined by a path all of whose edges
have odd labels.  Since $F_4$ is the only diagram with two connected
components with respect to this kind of connectivity, we deal with it
first.

\nl
$F_4$.  In view of Lemma 2.1 and symmetry of the diagram, we may
assume that $x_1$ has order greater than 2.  We describe a
calculational proof using GAP. The group $W$ has 25 conjugacy classes,
of which 17 contain elements $x_1$ of order at least $3$.  We first
compute a complete set of representatives of the conjugacy classes of
pairs $(x_1,x_2)$ of elements of $W$ such that $x_1$ and $x_2$ are of
order at least $3$, conjugate and satisfy $x_1x_2x_1=x_2x_1x_2$.  (For
each of the 17 representatives $x_1$, we compute the orbits of
$C_W(x_1)$ in the conjugacy class of $x_1$ of elements $x_2$ with
$x_1x_2x_1=x_2x_1x_2$, and take representatives.)  There are 53 such
representatives.  Since none of these generates $W$, we may assume
that $x_3$ and $x_4$ are nontrivial. For each pair $(x_1,x_2)$ we
compute representatives of the $C_W(x_1,x_2)$-orbits under conjugation
on $C_W(x_1)$ with $[x_2,x_3\rangle^4 =[x_3,x_2\rangle^4$.  There are
513 such triples (with $x_3\ne1$). Finally, for each such triple, we
compute representatives $x_4$ of each $C_W(x_1,x_2,x_3)$-orbit of
elements conjugate to $x_3$, contained in $C_W(x_1,x_2)$, and
satisfying $x_3x_4x_3 = x_4x_3x_4$.  There are $441$ such quadruples
$(x_1,x_2,x_3,x_4)$. Each of these generates a proper subgroup of $W$.
In particular, for $F_4$, there are no extraordinary epimorphisms.

\medskip 
For the remaining cases, all $x_i$ will be conjugate and hence have
sign $-1$.  Moreover, {from} the diagrams it will be clear that we may
assume $x$ to be injective when viewed as a map on $\{1,\ldots,n\}$
(observe that if $i$ and $j$ are nodes with $x_i=x_j$, then $x_k=x_i$
for each node adjacent to $i$ but not to $j$ [by an odd-labelled
edge], so by induction, $x$ would be a constant map, and so the
corresponding $\f$ would have a cyclic image, which does not occur).
By the same argument, we may even assume $x$ to be strong injective,
by which we mean that $\langle x_i\rangle \ne \langle x_j\rangle $ for
$i\ne j$.

For $J\subseteq \{1,\ldots,n\}$, denote by $X_J$ the set of all
$W$-conjugacy classes of strong injective maps $x : J\to W$, denoted by
$j\mapsto x_j$, such that the $x_j$ $(j\in J)$ are all from the same
conjugacy class in $W$, have order at least 3 and satisfy $ [x_i, x_j
\rangle^{m_{i,j}} = [x_j, x_i \rangle^{m_{i,j}} $
and $\sg(x_i) = -1$ for all $i,j\in
J$. The sets $X_J$ can be calculated in GAP by
similar computations to those described for $F_4$ above.  Since
$X_J=\emptyset$ for some $J\subseteq\{1,\ldots,n\}$, there are no
extraordinary epimorphisms. We give some of the statistics involved.

\nl $H_4$.
There are $34$ conjugacy classes, but 
$|X_{\{4\}}| = 8$.
We find
$|X_{\{2,4\}}| = 6$,
$|X_{\{1,2,4\}}| = 6$,
and
$|X_{\{1,2,3,4\}}| = 0$.

\nl $E_6$.
There are 25 conjugacy classes, but $|X_{\{1\}}|=8$.
Now $|X_{\{1,2,6\}}|=11$ and
$|X_{\{1,2,4,6\}}|=0$.

\nl $E_7$.
There are 60 conjugacy classes, but $|X_{\{1\}}|=25$,
$|X_{\{2,3,5,7\}}|=840$,
and $|X_{\{1,2,3,5,7\}}|=0$.

\nl $E_8$.
There are 112 conjugacy classes, but $|X_{\{1\}}|=43$,
$|X_{\{2,3,5,7\}}|=17700$,
and $|X_{\{2,3,5,7,8\}}|=0$.
$\ \square$

%\medskip
%With a little more effort, all morphisms $A\to W$ can be determined.
%For example, for $F_4$, there is...
%$$[x_1,x_2,x_3,x_4] = 
%[s_1s_4s_3s_2s_1s_3s_2s_3,
%s_2s_3s_2s_3s_4s_3s_2s_1,
%s_2s_3s_2s_4s_3s_2s_1s_3,
%s_1s_3s_4s_3s_2s_1s_3s_2]
%$$
\nl
{\bf Corollary 5.5.} {\it For each of the Coxeter graphs $H_4$, $F_4$,
$E_6$, $E_7$, $E_8$, the colored Artin group $CA$ is a characteristic
subgroup of $A$.} $\ \square$

%%%%%%%%%%%%%%%%%%%%%%%%%%%%%%%%%%%%%%%%%%%%%%%%%%%%%%%
\bigskip\nl
{\titre References}

\bigskip
\item{[Art]}
E. Artin,
{\it Braids and permutations},
Ann. of Math. {\bf 48} (1947), 643--649.

\smallskip
\item{[Bou]}
N. Bourbaki,
{\it Groupes et alg\`ebres de Lie, Chapitres 4, 5, et 6},
Hermann, Paris, 1968.

\smallskip
\item{[BS]}
E. Brieskorn, K. Saito,
{\it Artin-Gruppen und Coxeter-Gruppen},
Invent. Math. {\bf 17} (1972), 245--271.

\smallskip
\item{[Fra]}
W.N. Franzsen,
{\it Automorphisms of Coxeter groups},
Ph. D. Thesis, University of Sydney, 2001.

\smallskip
\item{[Lin]}
V.Y. Lin,
{\it Artin braids and related groups and spaces},
Algebra. Topology. Geometry, Vol. 17, pp. 159--227,
308, Akad. Nauk SSSR, Vsesoyuz. Inst. Nauchn. i Tekhn. Informatsii, Moscow, 1979.

\smallskip
\item{[Zin]}
V.M. Zinde,
{\it Some homomorphisms of the Artin groups of the series $B_n$ and $D_n$ into groups of the same
series and into symmetric groups},
Uspehi Mat. Nauk {\bf 32} (1977), 189--190.

%%%%%%%%%%%%%%%%%%%%%%%%%%%%%%%%%%%%%%%%%%%%%
\bigskip\nl
\line{
\vbox{\halign{#\hfill\cr
Arjeh M. Cohen\cr
Department of Mathematics and Computer Science\cr
Eindhoven University of Technology\cr
PO Box 513\cr
5600 MB Eindhoven\cr
THE NETHERLANDS\cr
\noalign{\smallskip}
{\tentt A.M.Cohen@tue.nl}\cr}}
\hskip2truecm
\vbox{\halign{#\hfill\cr
Luis Paris\cr
Laboratoire de Topologie\cr
Universit\'e de Bourgogne\cr
UMR 5584 du CNRS, BP 47870\cr
21078 Dijon cedex\cr
FRANCE\cr
\noalign{\smallskip}
{\tentt lparis@u-bourgogne.fr}\cr}}
\hfill}

\end